\documentclass[
fontsize=10pt,
]{scrartcl}
\addtokomafont{disposition}{\color{SaddleBrown}}

\usepackage[
final,
theoremdefs,
]{latexdev}
\usepackage{unixode}

\usepackage{authblk}

\title{Reductions of Operator Pencils}
\author{Olivier Verdier\thanks{\url{olivier.verdier@math.ntnu.no}}}
\affil{Department of Mathematical Sciences, NTNU, 7491 Trondheim, Norway}

\hypersetup{
pdfauthor={Olivier Verdier},
}

\usepackage{tikz}
\usetikzlibrary{matrix}
\tikzstyle{commdiag}=[matrix of math nodes, row sep=3em, column sep=2.5em, text height=1.5ex, text depth=0.25ex,ampersand replacement=\&]
\tikzstyle{exseq}=[commdiag, column sep=2em]

\tikzset{>=stealth}
\tikzstyle{normalcond}=[]

\newenvironment{pr}[1][\proofname]{
\begin{proof}[#1]%
}{\end{proof}
}

\usepackage{enumerate}
\newenvironment{thmenumerate}{\begin{enumerate}[{\normalfont (i)}]}{\end{enumerate}}

\renewcommand{\ker}{\operatorname{ker}}

\newcommand{\R}{\mathbf{R}}
\newcommand{\CC}{\mathbf{C}}
\newcommand{\NN}{\mathbf{N}}

\ifxetex
\renewcommand*\U{U}
\else
\newcommand*\U{U}
\fi
\newcommand{\W}{W}

\newcommand{\kE}{\ker\E}
\newcommand{\EU}{\overline{\E\U}}
\newcommand{\Ak}{\A\kE}

\newcommand{\AK}{\overline{\Ak}}
\newcommand{\EUAK}{\overline{\E\U + \Ak}}
\newcommand{\cored}{\cred{1}\ored{\ 1}}
\newcommand{\ocred}{\ored{1}\cred{\ 1}}

\newcommand*{\JU}{J_{\U}}
\newcommand*{\JW}{J_{\W}}

\newcommand*{\ored}[1]{^{#1}}
\newcommand*{\cred}[1]{_{#1}}

\newcommand*{\ps}[2]{\left(#1,#2\right)}

\newcommand{\dd}{\mathrm{d}}
\newcommand{\ddt}{\frac{\dd}{\dd t}}

\newcommand{\Do}{\mathsf{D}}
\newcommand{\A}{\mathsf{A}}
\newcommand{\E}{\mathsf{E}}
\newcommand{\Sop}{\mathsf{S}}

\newcommand*{\normH}[1]{\abs{1}}

\renewcommand{\H}{\mathrm{H}}
\newcommand*{\Leb}[1]{\mathrm{L}^{#1}}

\newcommand{\orth}{^{\perp}}
\newcommand{\Span}{\mathrm{span}}
\newcommand*{\coker}{\operatorname{coker}}

\newcommand{\ii}{\operatorname{i}}

\newcommand{\grad}{\operatorname{grad}}
\renewcommand{\div}{\operatorname{div}}

\newcommand{\Id}{\operatorname{Id}}

\newcommand{\inv}{^{-1}}

\renewcommand*{\Re}{\operatorname{Re}}

\newcommand*{\Uninemat}[3]{%
	\begin{tikzpicture}
		\matrix(m)[commdiag]
		{ \& 0 \& 0 \& 0 \& \\
		0 \& #1 \& 0\\
		0 \& #2 \\
		0 \& #3 \\
		 \& 0 \& 0 \& 0 \& \\};
		\path[->]
		(m-1-2) edge (m-2-2)
		(m-2-2) edge (m-3-2)
		(m-3-2) edge (m-4-2)
		(m-4-2) edge (m-5-2);
		\path[->]
		(m-1-3) edge (m-2-3)
		(m-2-3) edge (m-3-3)
		(m-3-3) edge (m-4-3)
		(m-4-3) edge (m-5-3);
		\path[->]
		(m-1-4) edge (m-2-4)
		(m-2-4) edge (m-3-4)
		(m-3-4) edge (m-4-4)
		(m-4-4) edge (m-5-4);
		\path[->]
		(m-2-1) edge (m-2-2)
		(m-2-2) edge (m-2-3)
		(m-2-3) edge node[auto] {$\A$} (m-2-4)
		(m-2-4) edge (m-2-5);
		\path[->]
		(m-3-1) edge (m-3-2)
		(m-3-2) edge (m-3-3)
		(m-3-3) edge node[auto] {$\A$} (m-3-4)
		(m-3-4) edge (m-3-5);
		\path[->]
		(m-4-1) edge (m-4-2)
		(m-4-2) edge (m-4-3)
		(m-4-3) edge node[auto] {$\A\cred{1}$}(m-4-4)
		(m-4-4) edge (m-4-5);
	\end{tikzpicture}
}

\newcommand*{\Wninemat}[4]{%
	\begin{tikzpicture}
		\matrix(m)[commdiag]
		{ \& #1 \& 0 \&\\
		0 \& #2 \& 0\\
		0 \& #3 \& 0 \\
		0 \& #4 \& 0 \\
		 \& 0 \& 0 \& 0 \&\\};
    \path[->]
    (m-1-2) edge (m-2-2)
    (m-2-2) edge node[auto] {$\A\ored{1}$} (m-3-2)
    (m-3-2) edge (m-4-2)
    (m-4-2) edge (m-5-2);
    \path[->]
    (m-1-3) edge (m-2-3)
    (m-2-3) edge node[auto] {$\A$} (m-3-3)
    (m-3-3) edge (m-4-3)
    (m-4-3) edge (m-5-3);
    \path[->]
    (m-1-4) edge (m-2-4)
    (m-2-4) edge node[auto] {$\A$} (m-3-4)
    (m-3-4) edge (m-4-4)
    (m-4-4) edge (m-5-4);
    \path[->]
    (m-2-1) edge (m-2-2)
    (m-2-2) edge (m-2-3)
    (m-2-3) edge (m-2-4)
    (m-2-4) edge (m-2-5);
    \path[->]
    (m-3-1) edge (m-3-2)
    (m-3-2) edge (m-3-3)
    (m-3-3) edge  (m-3-4)
    (m-3-4) edge (m-3-5);
    \path[->]
    (m-4-1) edge (m-4-2)
    (m-4-2) edge (m-4-3)
    (m-4-3) edge (m-4-4)
    (m-4-4) edge (m-4-5);
	\end{tikzpicture}
}

\begin{document}

\maketitle

\begin{abstract}
We study problems associated with an operator pencil, i.e., a pair of operators on Banach spaces.
Two natural problems to consider are linear constrained differential equations and the description of the generalized spectrum.
The main tool to tackle either of those problems is the reduction of the pencil.
There are two kinds of natural reduction operations associated to a pencil, which are conjugate to each other.
%One way is called control reduction, the second is called observation reduction.

Our main result is that those two kinds of reductions commute, under some mild assumptions that we investigate thoroughly.

Each reduction exhibits moreover a pivot operator.
The invertibility of all the pivot operators of all possible successive reductions corresponds to the notion of regular pencil in the finite dimensional case, and to the inf-sup condition for saddle point problems on Hilbert spaces.

Finally, we show how to use the reduction and the pivot operators to describe the generalized spectrum of the pencil.
\end{abstract}

\newcommand*\amsclassification{15A21, 15A22, 34A30, 47A10, 65L80}
\textbf{AMS Classification:} \amsclassification
% \subjclass[2010]{Primary: \amsclassification}

\tableofcontents

\section{Introduction}

Constrained differential equations, or differential algebraic equations, have been extensively studied, both from a theoretical and a numerical point view.
It is essential to study the linear case in order to apprehend the general case.
In that case, the object to be studied is a \emph{matrix pencil}, i.e., a pair of matrices.
In that context, the fundamental tool is that of \emph{Kronecker decomposition}.
%When the generalized spectrum of the pencil does not fill up the whole complex plane, the pencil is called regular.

The Kronecker decomposition theorem being arduous to prove, the concept of \emph{reduction} was gradually developed, first in \cite{Wong} for the study of regular pencils, then in \cite[\S~4]{Wilkinson} and \cite{vanDooren} to prove the Kronecker decomposition theorem.
It was later used in \cite{Verdier} for the study of other invariants.

It is also related to the \emph{geometric reduction} of nonlinear implicit differential equations as described in \cite{Reich} or \cite{RabierRheinboldt}. 
In the linear case, those coincide with the observation reduction, as shown in \cite{thesis}.
It is also equivalent to the algorithm of prolongation of ordinary differential equation in the formal theory of differential equations, as shown in \cite{Reid}.

In \cite{Wilkinson} and \cite{Verdier}, one considers also the \emph{conjugate} of the reduction, i.e., the operation obtained by transposing both matrices, performing a reduction and transposing again.
In order to make a distinction between both operations, we call the first one ``observation'' reduction, and the latter, ``control'' reduction.
The control reduction coincides with one step of the \emph{tractability chain}, as defined in \cite{Marz}.

Both reduction also appear in the context of ``linear relations''.
For a system of operators $(\E,\A)$ defined from $\U$ to $\W$, there are two corresponding linear relations, which are subspaces of $\U\times\U$ and $\W\times\W$.
These are respectively called the left and right linear relation (\cite[\S~6]{BaskakovSpectral}, \cite[\S~5.6]{BaskakovRepresentation}).
These linear relations correspond to the differential equations $\E u' + \A u = 0$ and $(\E u)' + \A u = 0$ respectively.
As one attempts to construct semigroup operators, it is natural to study the iterates of those linear relations.
That naturally leads to iterates of observation or control reduction.
There are two significant differences with what we do in this paper.
First, as the image of $\E$ is not necessarily closed, we will only consider its closure before pursuing the reduction.
This changes the notions of what is defined as reduced or not (see the examples in \autoref{sec_examples}).
Second, our main concern is the commutation of the two reduction procedure (\autoref{sec_commutativity}), not particularly the iterations of one type of reduction only.

The idea behind reduction operations is to produce a new, ``less implicit'' system from an implicit system.
When the pencil is not reducible anymore, it is equivalent to an ordinary differential equation.

%We are however not aware of any attempt to define the observation reduction on infinite dimensional systems.

In order to tackle linear constrained partial differential equations, we investigate the \emph{infinite dimensional case}, i.e., we replace the finite dimensional spaces by Banach spaces.
The Stokes and Maxwell equations may be naturally regarded as operator pencils.
Other examples of linear constrained partial differential equations include linearized elastodynamics in \cite{Simeon}, the Dirac equation in the nonrelativistic limit in \cite[\S~3]{ThallerCauchy}, linear PDEs as studied in \cite{Campbell}, \cite{Seiler}, \cite{Debrabant} and \cite[\S~6]{BaskakovSemigroup}.
There is a vast body of literature on that subject, and we refer to \cite{FaviniYagi}, \cite{BaskakovSemigroup}, \cite[\S~5]{BaskakovSpectral}, \cite{BaskakovRepresentation}, \cite[\S~4]{Tischendorf} and the references therein for more references.

We are interested in the kind of structures that may be preserved when the pencil is defined on Banach spaces.
For instance what can be said of the index, what is a regular pencil, what is the solvability of the corresponding abstract differential equation?

The reduction operations proves to be a very useful tool for this, as it may be defined on Banach spaces with virtually no modifications.
Each reduction naturally exhibits a \emph{pivot operator}, which is well defined in the Banach space case as well.

The invertibility of the pivot operators is essential in two different contexts.
First, in the finite dimensional case, the invertibility of all the pivot operators is exactly equivalent to the property of regularity of the pencil, as shown in \cite{Verdier}.
Second, for the Stokes equation, and all saddle point problems, the invertibility of the pivot operator is equivalent with the inf-sup condition.
In other word, \emph{a saddle point problem is a regular pencil if and only if the inf-sup condition is fulfilled}, which we show in \autoref{prop_infsup}.

In the finite dimensional case, regular pencils are pencils whose spectrum does not fill up the whole complex plane.
The corresponding property of regular pencils in Banach spaces is that the spectrum of the full pencil is equal to the spectrum of all the successively  reduced pencils, as we shall see in \autoref{thm_resolvent}.
If the pencil is not regular, then its spectrum is the whole complex plane.

Some attempts were made in \cite{Campbell} to define a notion of index for PDAE, but with no tangible conclusion or result.
The general attitude towards the index of an operator pencil is that one obtains a well defined index after spatial discretization.
The trouble with this approach is that the index thus obtained would generally depend on the choice of discretization.

We argue that there can in fact hardly exist any equivalent of the index in the infinite dimensional case.
The {index} of a regular finite dimensional pencil is defined from the Kronecker decomposition theorem.
Unfortunately, this decomposition is not available in the infinite dimensional case.
One observes however that for finite dimensional regular pencils, the number of observation reductions is the same as the number of control reductions, and that number gives a suitable definition of the index, which we could hope to extend to the infinite dimensional case.
As we shall see in \autoref{sec_mult_reg_ored}, this reasoning is not valid on Banach spaces.
In other words, the number of reductions of one type only is an unsatisfactory indicator of the structure of the system.
In particular, this means that the \emph{notion of index is not sufficient for regular operator pencils}.

As a result, to better describe the structure of an operator pencil, one has to describe the effect of the successive application of both kinds of reductions.
There is a possibly staggering amount of situations to consider, because the result of the successive application of reductions of different kinds generally leads to different pencils.

As it turns out, the situation is not that hopeless.

In the finite dimensional case, we observe indeed that the two kinds of reductions \emph{commute}.
We will show that under some general assumptions (studied in \autoref{sec_normality}), \emph{the two kinds of reductions also commute in the infinite dimensional case} (see \autoref{thm_commutativity}).
This property is essential because is considerably simplifies the description of all the possible reductions of the system.

\subsection{Outline of the Paper}

We start by defining the two possible reductions in \autoref{sec_reduction}.
We then show in \autoref{sec_findim} the relations with defects and Kronecker indices in the finite dimensional case.
What can be expected from the finite dimensional case is that the two types of reduction commute.
We proceed to show that this is indeed the case in \autoref{sec_commutativity}, under \emph{normality conditions}, that are studied in detail in \autoref{sec_normality}.

The rest of the paper is devoted to study examples and applications of reduction and of the Commutativity property.

In \autoref{sec_examples}, we study a multiplication operator system, and a saddle point problem. 
In particular, we will show that the inf-sup condition is none other than the invertibility of the pivot operator occuring in the reduction.

We then proceed to show how the pivot operators are related to the generalized resolvent set, thus making an analogy with the regular pencils in the finite dimensional case.

Finally, we study applications of reductions for linear problems in \autoref{sec_linprob}.

\subsection{Notations and Conventions}

\subsubsection{System}

The formal setting is the data of two Banach spaces $\U$ and $\W$, and two \emph{bounded operators} $\E$ and $\A$ having the same domain $\U$ and codomain $\W$.
\[
\E,\A \colon \U \longrightarrow \W
\]
Such a pair of operators, or operator pencil will be called a \emph{system} in the sequel.

\subsubsection{Cokernel}
At several occasions in the sequel we will need the definition of the \emph{cokernel} in the infinite dimensional case.

\begin{definition}
\label{def_cokernel}
We define the \emph{cokernel} of an operator $\E$ defined from $\U$ to $\W$ as
\[
\coker \E := \W/\EU
.
\]
\end{definition}

\subsubsection{Block Operator Notation}
\label{sec_block_not}

We will use an operator block notation to define operators from product of Banach spaces to product of Banach spaces.
For instance, if $\U = \U_1 \times \U_2$ and $\W = \W_1 \times \W_2$, and if the operators $A_{ij}$, are defined from $\U_j$ to $\W_i$ for $1 \leq  i,j \leq  2$, we define the operator
\[
	\A = \begin{bmatrix}
		A_{11} & A_{12}\\ A_{21} & A_{22}
	\end{bmatrix}
\]
as the operator
\[
	\U_1 \times \U_2 \ni (u_1,u_2) \longmapsto (A_{11} u_1 + A_{12} u_2, A_{21} u_1 + A_{22} u_2) \in \W_1 \times \W_2
	.
\]

\subsubsection{Equivalent Systems}

Another important concept is that of \emph{equivalence}.
Two systems are equivalent if a change of variables transforms one system into the other.
The precise definition of equivalence is the following.

\begin{definition}
	\label{def_eq_sys}
	Two systems $(\E,\A)$ with domain $\U$ and codomain $\W$,  and $(\overline\E,\overline\A)$ with domain $\overline{\U}$ and codomain $\overline{\W}$ are \emph{equivalent} if there exists invertible linear mappings $\JU$ from $\U$ to $\overline{\U}$ and $\JW$ from $\W$ to $\overline{\W}$ such that 
	$\E = \JW\inv\overline{\E}\JU$ and $\A = \JW\inv\overline{\A}\JU$.
\end{definition}

\section{Reduction}
\label{sec_reduction}

The main tool used in this article is the process of \emph{reduction} of a system $(\E,\A)$.
We proceed to define the two kinds of reductions, the \emph{observation reduction} and the \emph{control reduction}.

Moreover, for each type of reduction there corresponds a \emph{pivot operator}.
For saddle point problems, the invertibility of that operator is equivalent to the inf-sup condition (see \autoref{prop_infsup}).

\subsection{Operators defined by invariant subspaces}

We will use the following decomposition property of operators with respect to invariant subspaces.

\begin{proposition}
\label{prop_subops}
Consider an operator $\Sop$ defined from a Banach space $X$ to a Banach space $Y$.
Consider also a closed subspace $X'\subset X$ and a closed subspace $Y'\subset Y$ such that
\[
\Sop X' \subset Y'
.
\]
The operators $\Sop'$ and $[\Sop]$ are then uniquely defined by the requirement that the following diagram commutes.
\begin{center}
\begin{tikzpicture}
\matrix(m) [commdiag, column sep=3.5em, row sep=4em]
{0 \& X' \& X \& X/X' \& 0\\
 0 \& Y' \& Y \& Y/Y' \& 0 \\};
\path[->]
(m-1-2) edge (m-1-3)
(m-2-2) edge  (m-2-3)
(m-1-2) edge node[auto] {$\Sop'$} (m-2-2)
(m-1-3) edge node[auto] {$\Sop$} (m-2-3)
(m-1-3) edge  (m-1-4)
(m-2-3) edge (m-2-4)
(m-1-4) edge node[auto] {$[\Sop]$} (m-2-4)
(m-1-1) edge (m-1-2)
(m-2-1) edge (m-2-2)
(m-1-4) edge (m-1-5)
(m-2-4) edge (m-2-5);
\end{tikzpicture}
\end{center}
\end{proposition}

The proof is elementary.

\subsection{Observation Reduction}

If $\E$ is not surjective, then the equation $\ddt{}\E u + \A u = 0$ contains \emph{constraints}.
Intuitively, the ``amount'' of constraints is measured by the cokernel of $\E$, namely $\W/\EU$.
The variables satisfying the constraints are given by the space $\A\inv \EU$.
The idea of the observation reduction is to create a new system, without the original constraints, where the variable satisfy the constraints.
We are therefore led to consider a new system which is defined from $\A\inv\EU$ to $\EU$.
It is easy to see that it is possible.

The ``observation reduced'' system is therefore defined by the subspaces
\[
\W\ored{1} := \EU
\]
and
\[
\U\ored{1} := \A\inv \EU := \bigl\{\, u\in\U : \A u \in \EU \,\bigr \}
.
\]

Observe that we have $\A\U\ored{1} \subset \W\ored{1}$ and $\E\U\ored{1} \subset \W\ored{1}$, so we may use \autoref{prop_subops} to define new operators from $\U\ored{1}$ to $\W\ored{1}$ and from $\U/\U\ored{1}$ to $\W/\W\ored{1}$.
We thus define the reduced operators $\E\ored{1}$, $\A\ored{1}$ and the pivot operator $[\A\ored{1}]$ by the requirement that the following diagram commutes.
\begin{center}
\begin{tikzpicture}
\matrix(m) [commdiag, column sep=3.5em, row sep=4em]
{0 \& \U\ored{1} \& \U \& \U/\U\ored{1} \& 0\\
0 \& \W\ored{1} \& \W \& \W/\W\ored{1} \& 0 \\};
\path[->]
(m-1-2) edge  (m-1-3)
(m-2-2) edge  (m-2-3)
(m-1-2) edge node[auto] {$\E\ored{1},\A\ored{1}$} (m-2-2)
(m-1-3) edge node[auto] {$\E,\A$} (m-2-3)
(m-1-3) edge  (m-1-4)
(m-2-3) edge  (m-2-4)
(m-1-4) edge node[auto] {$0,[\A\ored{1}]$} (m-2-4)
(m-1-1) edge (m-1-2)
(m-2-1) edge (m-2-2)
(m-1-4) edge (m-1-5)
(m-2-4) edge (m-2-5);
\end{tikzpicture}
\end{center}

In other words, $\E\ored{1}$ and $\A\ored{1}$ are defined as restrictions of $\E$ and $\A$ on $\U\ored{1}$ to $\W\ored{1}$, and the pivot operator $[\A\ored{1}]$ is defined by the quotient of $\A$ from $\U/\U\ored{1}$ to $\W/\W\ored{1}$.
Clearly, the quotient of $\E$ defined from $\U/\U\ored{1}$ to $\W/\W\ored{1}$ is zero.

Note that this defines for any integer $k$ a system $(\E\ored{k},\A\ored{k})$, along with their domain $\U\ored{k}$ and codomain $\W\ored{k}$, with the convention that $(\E\ored{0},\A\ored{0}) := (\E,\A)$.

\subsection{Control Reduction}

The ``control'' reduction is conjugate to that of the observation reduction.
One considers variables which are not differentiated in the equation $\ddt{}\E u + \A u = 0$.
The space corresponding to those variables is $\kE$.
We interpret those variables as \emph{control variables}.
Intuitively, those control variable only have an influence on the space $\Ak$.
In fact, since we are considering Banach spaces, we have to consider instead the space $\AK$.
This little complication will lead to the normality assumption in \autoref{sec_normality}.

The idea of the control reduction is now to get rid of the spaces $\kE$ and $\AK$ by taking quotients.

This leads to the definition of the spaces
\[
\U\cred{1} := \U/\kE
,
\]
and 
\[
\W\cred{1} := \W/\AK
.
\]

Obviously we have $\E\kE \subset \AK$ and $\A\kE \subset\AK$, so we may use \autoref{prop_subops}.
The reduced operators $\E\cred{1}$, $\A\cred{1}$ and $[\A\cred{1}]$ are thus uniquely defined by the requirement that the following diagram commutes.

\begin{center}
\begin{tikzpicture}
\matrix(m) [commdiag, column sep=3.5em, row sep=4em]
{0 \& \kE \& \U \& \U\cred{1} \& 0\\
0 \& \AK \& \W \& \W\cred{1} \& 0 \\};
\path[->]
(m-1-2) edge  (m-1-3)
(m-2-2) edge  (m-2-3)
(m-1-2) edge node[auto] {$0,[\A\cred{1}]$} (m-2-2)
(m-1-3) edge node[auto] {$\E,\A$} (m-2-3)
(m-1-3) edge  (m-1-4)
(m-2-3) edge  (m-2-4)
(m-1-4) edge node[auto] {$\E\cred{1},\A\cred{1}$} (m-2-4)
(m-1-1) edge (m-1-2)
(m-2-1) edge (m-2-2)
(m-1-4) edge (m-1-5)
(m-2-4) edge (m-2-5);
\end{tikzpicture}
\end{center}

In other words, the pivot operator $[\A\cred{1}]$ is defined as the restriction of $\A$ from $\kE$ to $\AK$, and $\E\cred{1}$ and $\A\cred{1}$ are defined as quotient operators from $\U\cred{1}$ to $\W\cred{1}$.

Note that this defines for any integer $k\in\NN$ a system $(\E\cred{k},\A\cred{k})$, along with their domains $\U\cred{k}$ and $\W\cred{k}$, with the convention that $(\E\cred{0},\A\cred{0}) := (\E,\A)$.

\begin{remark}
	\label{rk_cont_lit}
	In the finite dimensional case, the control reduction is \emph{conjugate} to the observation reductions, with respect to the transposition.
	In other words, performing an observation reduction and transposing is the same as transposing and performing a control reduction.
	
\end{remark}

\subsection{Irreducible Systems}

Let us collect the definition and elementary properties of an irreducible system.

\begin{definition}
	A system $(\E,\A)$ is \emph{irreducible} if it is neither control-reducible nor observation-reducible.
\end{definition}

It is straightforward to characterize irreducible systems.

\begin{proposition}
\label{prop_irreducibility}
A system $(\E,\A)$ is irreducible if and only if  $\E$ is injective and has dense image, i.e., recalling \autoref{def_cokernel}, is such that
\[
\begin{split}
\kE = 0
, \\
\coker\E = 0
.
\end{split}
\]
\end{proposition}

\section{Finite Dimensional Case}
\label{sec_findim}

In this section we assume that $\U$ and $\W$ are \emph{finite dimensional}.
In this case, it is easy to describe the effect of the reductions defined in \autoref{sec_reduction}.

\subsection{Indices and Defects}

In the finite dimensional case, a system $(\E,\A)$ is characterised up to equivalence by its \emph{Kronecker indices} (see, e.g., \cite{Gantmacher}), so we may describe precisely the effect of both the observation and the control reduction.
In fact, it is easier to use equivalent invariants called the \emph{defects} (see \cite{Verdier})
Those defects are defined as follows.

We first define the \emph{constraint defect} $\alpha _1$.
That integer roughly represents the number of ``algebraic'' variables, i.e., constraint variables that are not involved at all in the differential equations.

In order to define the constraint defect $\alpha _1$ properly, we must define an auxiliary operator $[\E]$.
The operator $\E$ defines the operator $[\E]$ defined on the quotient space $\U/\U\ored{1}$ to the quotient space $\W\ored{1}/\W\ored{2}$ by the fact that the following diagram commutes.

\begin{center}
\begin{tikzpicture}
\matrix(m) [commdiag]
{\U \& \W\ored{1}\\
\U/\U\ored{1} \& \W\ored{1}/\W\ored{2} \\};
\path[->]
(m-1-1) edge node[auto] {$\E$} (m-1-2)
(m-2-1) edge node[auto] {$[\E]$} (m-2-2)
(m-1-1) edge (m-2-1)
(m-1-2) edge (m-2-2);
\end{tikzpicture}
\end{center}

Notice that the operator $[\E]$ is surjective, but not necessarily injective.
The constraint defect $\alpha _1$ is thus defined as
\[
	\alpha _1(\E,\A) := \dim \ker [\E]
.
\]

We now turn to the definitions of the \emph{control} and \emph{observation} defects.

The first \emph{observation defect} $\beta\ored{1}$ represents the number of ``empty equations'', i.e., equations of the form $0=0$.
More precisely, $\beta\ored{1}$ is defined as the dimension of the cokernel of the pivot operator $[\A\ored{1}]$, i.e.,
\begin{equation}
	\label{def_betaored}
	\beta\ored{1}(\E,\A) := \dim \coker [\A\ored{1}]
.
\end{equation}

The first \emph{control defect} $\beta\cred{1}$ represents the number of variables that are not present at all in the equations.
It is defined as the dimension of the kernel of the operator $[\A\cred{1}]$, i.e.,
\begin{equation}
	\label{def_betacred}
	\beta\cred{1}(\E,\A) := \dim\ker [\A\cred{1}]
.
\end{equation}

We then define iteratively 
\[\alpha _{k+1}(\E,\A) := \alpha _1(\E\ored{k},\A\ored{k})
	,
\]
\[
	\beta\ored{k+1}(\E,\A) := \beta\ored{1}(\E\ored{k},\A\ored{k})
	,
\]
and
\[\beta\cred{k+1}(\E,\A) := \beta\cred{1}(\E\cred{k},\A\cred{k})
	.
\]

When no reduction is possible, a system is irreducible.
Any system has a unique underlying irreducible system.
The equivalence class of that irreducible system with respect to invariance, together with all the defects $\alpha $, $\beta_*$ and $\beta^*$ completely determine a system.
Note that reductions keep the underlying irreducible system unchanged, so it suffices to describe the effect of reduction with the defects.

\subsection{Defects and Reductions}

According to our definition of the defects $\alpha $ and $\beta\ored{1}$, we immediately see that if a system has a constraint and observation defect sequence of $(\alpha _1,\alpha _2,\ldots)$ and $(\beta\ored{1},\beta\ored{2},\ldots)$, then the corresponding defect sequences of the observation-reduced system $(\E\ored{1},\A\ored{1})$ is just $(\alpha _2,\ldots)$, $(\beta\ored{2},\ldots)$.
It is more difficult to see however that the control defects $\beta\cred{*}$ are preserved by the observation reduction.

It is also possible to show that, symmetrically, the control reduction will shift the constraint defects $\alpha $ and the control defects $\beta\cred{*}$, but will preserve the observation defects $\beta\ored{*}$.

\subsection{Properties of Finite-Dimensional Pencils}

Let us gather some observations stemming from the fact that the defect completely determine a pencil.

\subsubsection{Commutativity of reductions}

The combination of two reduction of different kind should lead to equivalent systems.
In other words, the systems $(\E\cored,\A\cored)$ and $(\E\ocred,\A\ocred)$ have the same defects and the same underlying irreducible system, so they are equivalent. 
We will see in \autoref{sec_commutativity} that not only this is true in the infinite dimensional case, but also that the equivalence between those two systems is \emph{canonical}.
Canonical means that there are natural mappings that map $\U\cored$ and $\W\cored$ to $\U\ocred$ and $\W\ocred$ respectively.

To our knowledge, this is the first time that this fact is noticed, even in the finite dimensional case.

\subsubsection{Regular Pencils and Generalized Resolvent}
\label{sec_findim_reg}

The system $(\E,\A)$ is called a \emph{regular pencil} if its resolvent set is not empty (see \autoref{def_resolvent_set}).
It is shown in \cite[\S~3.8]{Verdier} that it is equivalent to all the $\beta$ defects being zero.

This is proved in two steps.

Consider the system $(\E,\A)$ and a reduced system $(\overline{\E},\overline{A})$, either via a control or an observation reduction.
One can show that the resolvent set of $(\E,\A)$ and $(\overline{\E},\overline{A})$ are either equal if the pivot operator is invertible, or the resolvent set of $(\E,\A)$ is empty.
This is still true in the infinite dimensional case, as we shall see in \autoref{thm_resolvent}.

In the finite dimensional case, all systems are reducible after a finite number of reductions.
It is straightforward to show that for an irreducible system, the resolvent set is not empty (because the spectrum of a matrix cannot fill the whole complex plane).

This last argument does not hold in the infinite dimensional case, for two reasons.
First, because a system need not be irreducible after a finite number of steps, second because the resolvent set of an irreducible system may be empty, as we show in \autoref{rk_empty_resolvent}.

\subsubsection{Regular Pencils and Equivalence of Reductions}

If the system $(\E,\A)$ is a \emph{regular pencil}, i.e., if all the $\beta$ coefficient vanish, then control and observation reductions lead to \emph{equivalent} systems. 
In particular, if a regular pencil is observation-irreducible, then it is also control-irreducible.
This is the meaning of the following proposition.

\begin{proposition}
\label{prop_findim_regnobeta}
Suppose that a \emph{finite dimensional} system $(\E,\A)$ is control-irreducible (i.e., $\U\cred{1} = \U$ and $\W\cred{1} = \W$) but not observation-irreducible.
Then one of the pivot operators $[\A\ored{k}]$ is not invertible.
\end{proposition}
\begin{pr}
If a system is control-irreducible, then in particular $\alpha _1 = 0$.
If it is not observation-reducible, it thus means that for some integer $k$ we have $\beta\ored{k} \neq  0$.
The claim is proved since $\beta\ored{k} = \dim \coker [\A\ored{k}]$.
\end{pr}

%There is however no reason to believe that those systems are canonically equivalent.

We will see in \autoref{sec_mult_reg_ored} that \autoref{prop_findim_regnobeta} \emph{cannot} be extended to the infinite dimensional case.

\section{Normality Assumption}
\label{sec_normality}

We will show in \autoref{sec_commutativity} that the two reductions defined in \autoref{sec_reduction} do commute, under some assumptions, summarized in \autoref{def_normality}.
In this section, we study those assumptions and the relation with the assumption that $\Ak$ is closed.

%% \begin{theorem}
%% The following conditions are equivalent.
%% \begin{thmenumerate}
%% \item $\A\kE + \EU = \AK + \EU$
%% \item the map $\A:\quad \kE \rightarrow \E/\EU$ has a closed image
%% \item $\A\kE/(\EU\cap\AK)$ is closed
%% \item $\A(\kE/\kE\ored{1})+\EU$ is a closed subspace of $\W/\E\U$
%% \end{thmenumerate}
%% If any of those conditions are fulfilled, we call the system \emph{normal}.
%% \end{theorem}

\subsection{Definition of Normality}

\begin{definition} \label{def_normality}
We will call a system \emph{normal} if the conditions
\begin{equation}
\label{eqnormsum}
\EU +  \Ak = \EUAK
\end{equation}
and
\begin{equation}
\label{eqnormint}
\overline{\EU \cap \A\kE} = \EU \cap \AK
\end{equation}
are fulfilled.
\end{definition}

We will also define the weaker condition
\begin{equation}
\EU + \AK = \EUAK \label{eqnormsumw}
.
\end{equation}

\begin{remark}
\label{rk_norm_inc}
The inclusions 
\[\EU +  \Ak \subset \EU + \AK \subset \EUAK
,
\]
and
\[
\overline{\EU \cap \A\kE } \subset \EU \cap \AK
\]
always hold.
In particular, notice that \eqref{eqnormsum} implies \eqref{eqnormsumw}.
\end{remark}

\begin{remark}
The normality assumptions are always fulfilled in the finite dimensional case.
\end{remark}

\begin{remark}
	\label{rk_normal_no_alpha}
	If $\Ak \subset \EU$ then the system $(\E,\A)$ is normal.
	In the finite dimensional case, this corresponds to the condition that 
	\[
		\alpha _1(\E,\A) = 0
	.
	\]
	The latter condition may be interpreted as the absence of ``algebraic variables'', i.e., the absence of variable that intervene solely in constraint equations.
\end{remark}

\subsection{Normality and Closedness of $\Ak$}

As we shall see, the space $\A\kE$ is often closed in applications (see \autoref{lma_AkerE}).
The following straightforward result is thus very useful to prove normality of a system.

\begin{proposition}
	\label{prop_AkE_norm}
	Consider a system $(\E,\A)$. If $\A\kE$ is closed, then $(\E,\A)$ is normal if and only if $\AK + \EU$ is closed.
\end{proposition}

Note also that the normality assumption is very close to the assumption on the closedness of $\A\kE$.
We study this assertion in details in the following remarks.

\begin{remark}
Let us observe that the closedness of $\Ak$ does not imply the normality assumption.
Suppose that a Banach space $\W$ has two closed subspaces $\mathcal{A}$ and $\mathcal{B}$ which sum is not closed.
% By taking the quotient of $\mathcal{A}\cap\mathcal{B}$ we may assume that $\mathcal{A}\cap\mathcal{B} = 0$.
Now consider the Banach space $\U := \mathcal{A} \times \mathcal{B}$.
We construct the system $(\E,\A)$ (using the block operator notation of \autoref{sec_block_not}) as
\[
	\E := \begin{bmatrix}
		\Id & 0
	\end{bmatrix}
	,
	\qquad
	\A := \begin{bmatrix}
		0 & \Id
	\end{bmatrix}
	.
\]
Clearly, the operators $\E$ and $\A$ are continuous from $\U$ to $\W$.

Now, observe that $\mathcal{B} = \A\ker\E$, which is closed, but the system is nevertheless not normal because $\Ak + \E\U$ is not closed.
\end{remark}

\begin{remark}
In a similar fashion, we observe that the normality assumptions do not imply the closedness of $\Ak$.
Consider a Banach space $\U\ored{1}$ continuously and densely injected in a Banach space $\W$ via an injection mapping $\ii$.
Build now $\U := \U\ored{1} \times \U\ored{1}$, and define 
\[
	\E := \begin{bmatrix}
		\ii & 0
	\end{bmatrix}
	,
	\qquad
	\A := \begin{bmatrix}
		0 & \ii
	\end{bmatrix}
	.
\]
It is clear that $\Ak = \E\U$, so the normality assumptions are fulfilled, but $\Ak$ is not closed.

Note however that the normality assumptions almost imply the closedness of $\Ak$ as the following proposition shows.

\begin{proposition}
Consider the condition
\begin{equation}\label{eqnormints}
\A\kE \cap \EU = \AK \cap \EU
,
\end{equation}
which is stronger that \eqref{eqnormint}.
Assume that a system $(\E,\A)$ fulfills both \eqref{eqnormsum} and \eqref{eqnormints}.
The subspace $\A\kE$ must then be closed.
\end{proposition}
\begin{pr}
Indeed, if $z\in\AK$, then, using \eqref{eqnormsum}, it may be written as
\[
z = \A k + y
,
\]
where $k\in\kE$ and $y\in\EU$.
Now since $y\in\EU\cap\AK$, we may use \eqref{eqnormints} to obtain that $y$ may be written $y=\A k_0$, where $k_0\in \kE$.
As a result we obtain $z=\A(k+k_0)$, so $\Ak$ is closed.
\end{pr}
\end{remark}

\subsection{Equivalent conditions}

Let us study some equivalent formulations of the normality assumptions of \autoref{def_normality}.

\subsubsection{First Normality Condition}

We are going to show in \autoref{prop_normal_seq} that the condition \eqref{eqnormsum} is connected to the exactness of a sequence which is always exact in the finite dimensional case.

Before proceeding further, we need a more useful descriptions of $\ker\E\ored{1}$ and $\coker\E\cred{1}$.
\begin{lemma}\label{lma_cokerEc}
For a system $(\E,\A)$, we have
\begin{itemize}
	\item $\ker\E\ored{1} = \ker\E \cap \U\ored{1}$
	\item $\coker\E\cred{1} \equiv \W/(\overline{\E\U + \Ak})$
\end{itemize}
\end{lemma}
\begin{pr}
The first observation is obvious. For the second one, notice that 
$y+\AK \in \W\cred{1}$ if and only if there exists a sequence $x_n\in\U$ such that
\[\lim_{n\rightarrow \infty } \inf_{z\in\AK} \|y - \E x_n - z \| = 0 
.
\]
It implies that there exists a sequence $k_n\in\kE$ such that $\E x_n + \A k_n$ converges towards $y$ so $y\in\overline{\E\U + \Ak}$.
On the other hand, if $y\in\overline{\E\U + \Ak}$ then there exists sequences $x_n\in\U$ and $k_n\in\kE$ such that $\E x_n + \A k_n$ converges towards $y$, which implies that
\[
\inf_{z\in\Ak} \| y - \E x_n - z\| \leq  \| y - \E x_n - \A k_n \| \rightarrow 0
,
\]
so the claim is proved.
\end{pr}

\begin{proposition}\label{prop_normal_seq}
Condition \eqref{eqnormsum} is equivalent to the exactness of the following sequence.
\begin{center}
\begin{tikzpicture}[ampersand replacement=\&]
\matrix(m) [commdiag,]
{0 \& \ker\E\ored{1} \& \kE \& \coker\E \& \coker\E\cred{1} \& 0\\};
\path[->]
(m-1-1) edge (m-1-2)
(m-1-2) edge (m-1-3)
(m-1-3) edge node[auto] {$\A$} (m-1-4)
(m-1-4) edge (m-1-5)
(m-1-5) edge (m-1-6);
\end{tikzpicture}
%% \begin{tikzpicture}[ampersand replacement=\&]
%% \matrix(m) [commdiag,]
%% {0 \& \ker\E\ored{1} \& \kE \& \AK/(\AK\cap\EU) \& 0\\};
%% \path[->]
%% (m-1-1) edge (m-1-2)
%% (m-1-2) edge (m-1-3)
%% (m-1-3) edge node[auto] {$\A$} (m-1-4)
%% (m-1-4) edge (m-1-5);
%% \end{tikzpicture}
%% \begin{tikzpicture}[ampersand replacement=\&]
%% \matrix(m) [commdiag,]
%% {0 \& \ker\E/\ker\E\ored{1} \& \coker\E \& \coker\E\cred{1} \& 0\\};
%% \path[->]
%% (m-1-1) edge (m-1-2)
%% (m-1-2) edge node[auto] {$\A$} (m-1-3)
%% (m-1-3) edge  (m-1-4)
%% (m-1-4) edge (m-1-5);
%% \end{tikzpicture}
\end{center}
\end{proposition}

\begin{pr}
The exactness of the sequence is clear except for the fact that $\A$ maps $\kE$ onto $\coker\E\cred{1}$.
Using \autoref{lma_cokerEc}, this is equivalent to $\overline{\E\U + \Ak} \subset \EU + \Ak$, which, considering \autoref{rk_norm_inc}, is exactly \eqref{eqnormsum}.
%% The exactness of the second sequence is precisely equivalent to
%% \[
%% \overline{\E\U + \Ak} \subset \EU + \Ak
%% .
%% \]
%% The exactness of the first sequence implies
%% \[
%% \AK \cap \EU \subset \A \ker\E\ored{1}
%% ,
%% \]
%% and we conclude with \autoref{lma_AkEU}.
%% 
%% For the other direction, we now only need to prove that $\A$ is surjective on $\AK/(\AK\cap\EU)$.
%% Take $z \in\AK$.
%% We may write $z = w + \A k$, where $w \in\EU$.
%% This implie in turn that $w\in\AK$, so $w\in\AK\cap\EU$ and thus $z + \AK\cap\EU = \A k + \AK\cap\EU$.
\end{pr}

\subsubsection{Second Normality Condition}

In order to give an equivalent formulation of condition \eqref{eqnormint}, we need to describe the space $\Ak\cap\EU$.

%We will need the following Lemma anyway during the course of the proof of \autoref{lma_pivot_ocred}.
\begin{lemma}\label{lma_AkEU}
For a system $(\E,\A)$, we have
\[
\A\ker \E\ored{1} = \Ak \cap \EU
.
\]
\end{lemma}
\begin{pr}
If $y \in \A\ker\E\ored{1}$ then there exists $k\in\kE$ such that $y=\A k$ and $\A k \in\EU$, so clearly $y \in \Ak \cap \EU$.
On the other hand if $y\in\Ak \cap \EU$ then $y = \A k$ for $k\in\kE$ and $y \in \EU$, hence $k\in\U\ored{1}$ from which we get $k\in\ker\E\ored{1}$ and the result is proved.
\end{pr}

It is now clear that, in view of \autoref{lma_AkEU}, the assumption \eqref{eqnormint} is equivalent to
\[
\EU \cap \AK \subset \overline{\A\kE\ored{1}}
.
\]

\section{Commutativity of the Reductions}
\label{sec_commutativity}

We turn to the most important section of this paper
and set out to prove that the systems $(\E\cored,\A\cored)$ and $(\E\ocred,\A\ocred)$ are \emph{canonically} equivalent.

\subsection{Natural mappings $\JU$ and $\JW$}

In order to establish the equivalence, we define two natural maps.
The first one, $\JU$, maps $\U\ocred$ to $\U\cred{1}$.
The second map, $\JW$, maps $\W\ocred$ to $\W\cred{1}$.

We define the operator $\JU$ as follows.
Pick $u + {\ker\E\ored{1}} \in \U\ocred$, where, by definition of $\U\ocred$, $u\in \U\ored{1}$.
Now since $\kE\ored{1} \subset \kE$, we may map $u + {\ker\E\ored{1}}$ into $\U\cred{1}$.
Thus we obtain a linear, continuous mapping
\[
\JU \colon  \U\ocred \rightarrow \U\cred{1}
,
\]
such that $\|\JU\| \leq  1$.

Similarly, we define the operator $\JW$ as follows.
Pick an element $w+\overline{\A\ker\E\ored{1}} \in \W\ocred$.
Since $\overline{\A\ker\E\ored{1}} \subset \AK$, we may map $w + \overline{\A\ker\E\ored{1}}$ to $w + \AK \in \W\cred{1}$.
Thus we obtain a linear, continuous mapping
\[
\JW\colon  \W\ocred \rightarrow \W\cred{1}
,
\]
such that $\|\JW\| \leq  1$.

\begin{remark}
	The norms of $\JU$ and $\JW$ are \emph{intrinsic} properties of the system.
In the finite dimensional case, those numbers would be partial indicators of how well conditioned the transformation to the Kronecker canonical form is.
\end{remark}

We now show that under the assumptions studied in \autoref{sec_normality}, the mappings $\JU$ and $\JW$ are isomorphisms to $\U\cored$ and $\W\cored$ respectively.

\subsection{Commutativity}

We break down the proof in several lemmas, each showing exactly which assumptions are necessary for what.

\subsubsection{Reduced Subspaces}

%The following observation is proved in a similar way.

%\begin{lemma}\label{prop_Wcored}
%If $y + \AK \in \W\cred{1}$, then 
%\[
%y + \AK \in \W\cored \iff y \in \EUAK
%.
%\]
%\end{lemma}
%\begin{pr}
%The proof goes along the same line as the proof of \autoref{lma_cokerEc}.
%\end{pr}

\begin{lemma}
\label{lma_JU}
The mapping $\JU$ is injective and its image is included in $\U\cored$, i.e.,
\[
\JU(\U\ocred) \subset \U\cored
.
\]
Moreover, if  \eqref{eqnormsum} is fulfilled, its image is exactly $\U\cored$, i.e.,
\[
\JU(\U\ocred) = \U\cored
.
\]
\end{lemma}

\begin{pr}
%% We build a map $\JU$ from $\U\ocred$ to $\U\cred{1}$, and we will show that $\JU$ is injective, that $\JU(\U\ocred) = \U\cored$ and that it is an isometry.
\begin{enumerate}
\item $\JU$ is injective, because if $\JU(u) = 0$, then $u\in\kE$, but since $u\in\U\ored{1}$, we conclude that $u \in\kE\cap \U\ored{1} = \ker\E\ored{1}$.
\item If $u\in\U\ored{1}$ then $\A\cred{1}(u+\kE) \subset  \overline{\E \U} + \overline{\A \kE} \subset \overline{\E\U + \Ak}$, so $u + \ker{\E} \in \U\cred{1}$.
We conclude that $\JU(\U\cored) \subset \U\cred{1}$.
\item Using \eqref{eqnormsum},
if $\A(u) \in \overline{\E\U + \Ak} \subset \EU + \Ak$, then there exists $k\in\kE$ such that $\A(u-k)\in\EU$ so $u-k \in \U\ored{1}$, so $\JU(\U\cored) = \U\cred{1}$.
\end{enumerate}
\end{pr}

\begin{lemma}
	\label{lma_JW}
The mapping $\JW$ maps $\W\ocred$ into $\W\cored$, i.e.,
\[
\JW(\W\ocred) \subset \W\cored
.
\]
Moreover, if \eqref{eqnormsumw} is fulfilled then 
\[\JW(\W\ocred) = \W\cored
.
\]
If \eqref{eqnormint} is fulfilled, then $\JW$ is injective.
\end{lemma}
\begin{pr}
\begin{enumerate}
\item Take $y+\overline{\A\kE\ored{1}} \in \W\ocred$. It is mapped to $y+\AK$, since $\overline{\A\kE\ored{1}}\subset \AK$.
\item Suppose that $\JW(y+\overline{\A\kE\ored{1}})=0$, i.e., $y\in\AK$. 
Assumption \eqref{eqnormint} allows to conclude that $y\in\overline{\A\kE\ored{1}}$ so $\JW$ is injective.
\item The image by $\JW$ of an element $y\in\W\ored{1}$ belongs to $\EU + \AK \subset \overline{\E\U + \Ak}$ hence $\JW(\W\ocred) \subset \W\ored{1}$.
\item Suppose that $y\in\W\cored$. This means that $y\in \overline{\E\U + \Ak}$.
Using \eqref{eqnormsumw} we obtain $y\in \EU + \AK$, so $y + \AK \in \EU$ and we conclude that $\JW(\W\ocred) = \W\cored$.
\end{enumerate}
\end{pr}

\subsubsection{Pivot Operators}

\begin{lemma}
\label{lma_pivot_cored}
We have
\[
[\A\ored{1}]\quad\text{invertible} \implies [\A\cored]\quad\text{invertible}
.
\]
Moreover, if \eqref{eqnormsum} is fulfilled, then
\[
[\A\cored]\quad\text{invertible} \implies [\A\ored{1}]\quad\text{invertible}
.
\]
\end{lemma}

\begin{pr}
	Notice first that $[\A\ored{1}]$ is injective.
%Let us show that $[\A\ored{1}]$ invertible is equivalent to $[\A\cored]$ invertible.
\begin{enumerate}
	\item Consider $y+\EUAK \in \coker\E\cred{1}$; since $[\A\ored{1}]$ is invertible, there exists $x\in\U$ such that $\A x - y \in \EU$, so $\A x + \EU + \AK = y + \EU + \AK$.
	\item
	Consider $y+\EU$.
	By projecting on $\coker\E\cred{1}$ and using that $[\A\cored]$ is invertible, we obtain $x\in\U$ such that $\A x = y + \EUAK$.
Using the assumption \eqref{eqnormsum}, there exists $k\in\kE$ such that $\A(x+k) = y + \EU$.
\end{enumerate}
\end{pr}

\begin{lemma}
\label{lma_pivot_ocred}
We have $\ker[\A\cred{1}] = \ker[\A\ocred]$, which implies in particular
\[
[\A\cred{1}]\quad\text{injective} \iff [\A\ocred]\quad\text{injective}
.
\]
Moreover,
\[
[\A\cred{1}]\quad\text{surjective} \implies [\A\ocred]\quad\text{surjective}
.
\]
If \eqref{eqnormsum} and \eqref{eqnormint} are fulfilled, then
\[
[\A\ocred]\quad\text{surjective} \implies [\A\cred{1}]\quad\text{surjective}
.
\]
\end{lemma}
\begin{pr}
	First we show that $\ker[\A\cred{1}] = \ker[\A\ocred]$.
\begin{enumerate}
%% \item If $[\A\ocred]$ is not injective, then there exists $k\in\kE$ such that $\A k = 0$ so $\A\cred{1}$ is not injective.
\item 
	Pick an element $x\in\ker[\A\ocred]$.
	Then $x\in\ker\E\ored{1}$ is such that $\A x = 0$, so $x\in\ker[\A\cred{1}]$.
\item 
	Pick $x\in\ker \A\cred{1}$. 
	Then $x\in\U\ored{1}$, so $x\in\kE\ored{1}$ and $x\in\ker[\A\ocred]$.
\end{enumerate}
Next we show that $[\A\cred{1}]$ has a closed image if and only if $[\A\ocred]$ has a closed image.
\begin{enumerate}
\item Assume that $[\A\cred{1}]$ has a closed image, i.e., that $\A\kE$ is closed.
Using \autoref{lma_AkEU}, we obtain that $\A\ker\E\ored{1}$ is closed.
\item Assume that $[\A\ocred]$ has a closed image, i.e., that $\A\kE\ored{1}$ is closed.
Take $y\in\overline{\A\kE}$.
Using \eqref{eqnormsum}, there exists $k\in\kE$ such that $y-\A k \in \EU$.
As a result, $y-\A k \in \EU \cap \AK$, so
using \eqref{eqnormint}, $y-\A k \in \overline{\A\kE\ored{1}}$, so this means that $y-\A k = \A k_0$ with $k_0\in\kE\cap\U\ored{1}$, so the image of $[\A\cred{1}]$ is closed.
\end{enumerate}
\end{pr}

\subsubsection{Commutativity Theorem}

We collect the result of the preceding Lemmas.

\begin{theorem}
\label{thm_commutativity}
If a system $(\E,\A)$ is normal (\autoref{def_normality}), then the systems $(\E\cored,\A\cored)$ and $(\E\ocred,\A\ocred)$ are equivalent (\autoref{def_eq_sys}).
Moreover, we have
\begin{equation}
	\label{eq_Aored_Acored}
[\A\ored{1}]\quad\text{invertible}\iff [\A\cored]\quad\text{invertible}
,
\end{equation}
and
\begin{equation}
	\label{eq_Acred_Aocred}
[\A\cred{1}]\quad\text{invertible} \iff [\A\ocred]\quad\text{invertible}
.
\end{equation}
\end{theorem}
\begin{pr}
	Under the normality assumptions, we may apply \autoref{lma_JU} and \autoref{lma_JW}.
	As a result, $\JU$ and $\JW$ are continuous invertible operators, and by the Banach theorem, their inverse is continuous as well.
	This shows that $\E\ocred = \JW\inv\E\cored\JU$ and that $\A\ocred = \JW\inv\A\cored\JU$, so the systems $(\E\ocred,\A\ocred)$ and $(\E\cored,\A\cored)$ are equivalent.
	The equivalences \eqref{eq_Aored_Acored} and \eqref{eq_Acred_Aocred} are consequences of \autoref{lma_pivot_ocred} and \autoref{lma_pivot_cored}.
\end{pr}

\subsection{Exact Sequences}

It is fruitful to redefine the reductions using the language of exact sequences.
This gives a feeling of the reasons behind the commutativity of the reductions.

%We will follow the following conventions for the sake of conciseness of notations.
In the following diagrams, the arrows with no labels are natural injections (from a subspace to an ambient space), or natural projections (from a space to a quotient space).
The arrows labeled ``$\A$'' are combinations of $\A$ with either a natural injection or a natural projection.

With these conventions in mind we may \emph{define} $\U\ored{1}$ and $\W\ored{1}$ by requiring the exactness of the two following diagrams.
\begin{center}
\begin{tikzpicture}
\matrix(m) [commdiag,]
{0 \& \U\ored{1} \& \U \& \coker \E \& {\coker[\A\ored{1}]} \& 0\\};
\path[->]
(m-1-1) edge (m-1-2)
(m-1-2) edge (m-1-3)
(m-1-3) edge node[auto] {$\A$} (m-1-4)
(m-1-4) edge (m-1-5)
(m-1-5) edge (m-1-6);
\end{tikzpicture}

\begin{tikzpicture}
\matrix(m) [commdiag,]
{0 \& \W\ored{1} \& \W \& \coker \E \& 0\\};
\path[->]
(m-1-1) edge (m-1-2)
(m-1-2) edge (m-1-3)
(m-1-3) edge  (m-1-4)
(m-1-4) edge (m-1-5);
\end{tikzpicture}
\end{center}

If $\Ak$ is closed, we may similarly define the spaces $\U\cred{1}$ and $\W\cred{1}$ by the exactness of the following diagrams.

\begin{center}
\begin{tikzpicture}
\matrix(m) [commdiag, , row sep=4em]
{0 \& \kE \& \U \& \U\cred{1} \& 0\\};
\path[->]
(m-1-1) edge (m-1-2)
(m-1-2) edge (m-1-3)
(m-1-3) edge  (m-1-4)
(m-1-4) edge (m-1-5);
\end{tikzpicture}

\begin{tikzpicture}
\matrix(m) [commdiag,]
{0 \& {\ker [\A\cred{1}]} \& \kE \& \W \& \W\cred{1} \& 0\\};
\path[->]
(m-1-1) edge (m-1-2)
(m-1-2) edge (m-1-3)
(m-1-3) edge[normalcond] node[auto] {$\A$} (m-1-4)
(m-1-4) edge  (m-1-5)
(m-1-5) edge (m-1-6);
\end{tikzpicture}
\end{center}

Recall that in the finite dimensional case, the spaces $\ker[\A\cred{1}]$ and $\coker[\A\ored{1}]$, appearing in the sequences above, are related to the defects by \eqref{def_betaored} and \eqref{def_betacred}.

For normal systems, those exact diagrams may be interwoven with each other.
This is illustrated in \autoref{fig_commutativity}.

\begin{figure}[H]
	\centering
\Uninemat{\kE\ored{1} \& \kE \& \AK/(\AK\cap\EU)}%
{\U\ored{1} \& \U \& \coker \E \& \coker[\A\ored{1}] }%
{\U\ored{1}\cred{1} \& \U\cred{1} \& \coker{\E\cred{1}} \& \coker[\A\cored]}

\Wninemat{\ker[\A\ocred] \& \ker[\A\cred{1}]}%
{\kE\ored{1} \& \kE \& \kE/\kE\ored{1}}%
{\W\ored{1} \& \W \& \coker\E}%
{\W\ored{1}\cred{1} \& \W\cred{1} \& \coker\E\cred{1}}
\caption{
This figure is a summary of the results of \autoref{sec_commutativity}, when the normality assumption holds.
In the infinite dimensional case, the rows and columns need not be exact.
In the \emph{finite} dimensional case, however, all the rows and columns are exact, without any extra assumption, and the commutativity of the reduction may be proved by simple diagram chasing.
}
\label{fig_commutativity}
\end{figure}

\section{Examples}
\label{sec_examples}

Let us begin the study of some examples by a remark that will help to compute the quotient spaces involved.

\begin{remark}
In Hilbert spaces, there is an easy way to compute the quotient by a subspace.
Say that $\H$ is a Hilbert space that admits the topological decomposition
\[
\H =  A \oplus B
,
\]
where $A$ and $B$ are two closed subspaces of $\H$.
Then one can check that
\[
\H/A \equiv B
.
\]
This is because on the one hand, $\H/A$ is isometric to $A\orth$, and on the other hand $A\orth$ is isomorphic to $B$.
\end{remark}

\subsection{Multiplication Operator Example}
\label{sec_multiplication_operator}

We consider the case where $\E$ is a multiplication operator by the characteristic function of an interval.
The functional spaces are
\begin{equation*}
\U := \H^1 := \H^1(\R),
\qquad 
\W := \Leb2 := \Leb2(\R)
,
\end{equation*}
and the operators are
\begin{equation*}
\E := m(x) \ii,
\qquad 
\A = D_x
.
\end{equation*}
The operator $\ii$ is the injection of $\H^1(\R)$ into $\Leb2(\R)$, and the function $m$ is defined by 
\[
	m(x) = 1 - \chi _{J}
	,
\]
where $\chi _{J}$ is the characteristic function of a finite interval 
\[
J := (a,b)
.
\]

So the action of $\E$ on a function $u\in\H^1$ is simply given by
\[
\E u  = (1-\chi _J) u
.
\]

This example was studied in \cite[Example~2.1]{FaviniYagi}.

We compute
\[
\E\U = \{\, u\in\H^1 : u_{|J} = 0 \,\}
,
\]
so we obtain readily
\begin{equation}
	\label{eq_mult_U1W1}
\W\ored{1} := \EU = \{\, u\in\Leb2 :  u|_{J} = 0 \,\}
,
\qquad
\U\ored{1} = \{\, u\in\H^1 : u_{x}|_{J} = 0 \,\}
.
\end{equation}

$\E\ored{1}$ is the restriction of $\E$ on $\U\ored{1}$, with codomain $\W\ored{1}$.
We see that $\E\ored{1}$ is injective and has dense image.

In the sequel we will make the identification
\[
\W\ored{1} \equiv \Leb2(\neg J)
,
\]
where $\neg J$ is the set
\[
\neg J := (-\infty ,a) \cup (b,+\infty )
.
\]

% Assume that $J = (a,b)$, finite interval.

Let us define the continuous linear form $\varphi $ defined on $\H^1$ by
\[
\varphi (u) := u(a) - u(b)
.
\]

We see that $\ker \varphi $ is naturally decomposed into
\[
\ker \varphi  = \U\ored{1} \oplus \H^1_0(J)
.
\]

By choosing a supplementary space to $\ker \varphi $ in $\U$, we may decompose $\U$ as
\begin{equation}
\label{eq_decomposition_U}
\U = \U\ored{1} \oplus \H_0^1(J) \oplus \R
,
\end{equation}
where $\R$ denotes a one-dimensional subspace complementary to $\ker \varphi $, for example $(\ker \varphi )\orth$, or the span of any function $u\in\H^1$ such that $u(a)\neq u(b)$.
%where $\H_a(J)$ is defined as
%[
%\H_a(J) := \bigl\{ u\in\H^1(J):\ u(a) = 0\bigr\}
%.
%\]

So we obtain
\[
\U/\U\ored{1} \equiv \H^1_0(J) \oplus \R
.
\]

We also readily see that
\[
\W/\W\ored{1} \equiv \Leb2(J)
.
\]

%The operator $[\A\ored{1}]$ is just $\A$ restricted on $\H_a$, so $[\A\ored{1}]$ is invertible, and inf-sup condition is satisfied.
Let us compute further
\[ \ker\E = \bigl\{\, u\in\H^1 : u|_{\neg J} = 0 \,\bigr\} \equiv \H^1_0(J)\]
and
\[
\A\ker\E \equiv \Leb2_0(J) := \Bigl\{\, u\in\Leb2(J) : \int_J u = 0 \,\Bigl\}
.
\]
%otherwise $\A\ker\E \equiv \Leb2(J)$.
The set $\Ak$ is thus closed in $\Leb2$. Incidentally, this shows that the corresponding pivot operator $[\A\cred{1}]$ is invertible, since $\A$ is injective on $\ker\E$.

Moreover since $\Ak$ is closed, using the observation of \autoref{prop_AkE_norm}, the normality assumptions reduces to the closedness of $\EU + \Ak$, which is straightforward.
We conclude that the system $(\E,\A)$ is normal.

Using \eqref{eq_decomposition_U}, i.e., $\U = \U\ored{1} \oplus \ker\E \oplus \R$, we obtain
\begin{equation}\label{eq_mult_cred}
\U\cred{1} := \U/\ker\E \equiv \R \oplus \U\ored{1}
,
\qquad
\W\cred{1} = \Leb2/\Leb2_0(J) \equiv \R \oplus \Leb2(\neg J)
,
\end{equation}
where $\R$ denotes here the one-dimensional subspace of $\Leb2(\R)$ spanned by the characteristic function $\chi _J$ of the interval $J$.

$\E\cred{1}$ is injective, so no control reduction is possible.
We see however that the image of $\E\cred{1}$ is not dense and that $\W\cred{1}/\overline{\E\cred{1}\U\cred{1}} \equiv \R$.
We may thus perform an observation reduction, which leads to the spaces
\[
\U\cored = \U\ored{1}
\]
and
\[
\W\cored = \W\ored{1}
.
\]

The diagrams of \autoref{fig_commutativity} in the multiplication operator case are represented on \autoref{fig_multop}.

\begin{figure}[H]
	\centering
\Uninemat{0 \& \H^1_0(J) \& \Leb2_0(J)}{\U\ored{1} \& \H^1 \& \Leb2(J) \& 0}{\U\ored{1} \& \R \oplus \U\ored{1} \& \R \& 0}
\Wninemat{0 \& 0}{0 \& \H^1_0(J) \& \H^1_0(J)}{\Leb2(\neg J) \& \Leb2 \& \Leb2(J)}{\Leb2(\neg J) \& \R \oplus \Leb2(\neg J) \& \R}
\caption{
	The diagrams of \autoref{fig_commutativity} in the case of the multiplication operator example of \autoref{sec_multiplication_operator}.
	All the rows and columns are exact.
}
\label{fig_multop}
\end{figure}

\subsection{The Reduced Multiplication Operator System}
\label{sec_mult_reg_ored}

Here we study in more details the ``control-reduced'' system obtained in \autoref{sec_multiplication_operator}.
The system we are studying is thus the one defined by \eqref{eq_mult_cred}, i.e., the system $(\E,\A)$ we are considering in this section is the one denoted $(\E\cred{1},\A\cred{1})$ in \autoref{sec_multiplication_operator}.

We proceed to show that the conclusions of \autoref{prop_findim_regnobeta} do not hold for that system.

Following \eqref{eq_mult_cred}, the spaces $\U$ and $\W$ are now defined as
\[
	\U := \bigl\{\, u\in \H^1 : u_{xx} = 0\quad\text{on $J$} \,\bigr\}
	,
\]
\[
	\W := \bigl\{\, u\in\Leb2 : u_x = 0\quad\text{on $J$} \,\bigr\}
	.
\]
The operators $\E$ and $\A$ are defined as in \autoref{sec_multiplication_operator}.
The operator $\E$ is injective, but its image is not dense.
Observe that $\U\ored{1}$ and $\W\ored{1}$ are still defined as in \eqref{eq_mult_U1W1}.
With the notations of \autoref{sec_multiplication_operator}, we have
\[
	\U/\U\ored{1} \equiv  (\ker \varphi )\orth \equiv \R
	,
\]
\[
	\W/\W\ored{1} \equiv \Span (\chi _J) \equiv \R
	,
\]
and the pivot operator $[\A\ored{1}]$ is a non-zero operator from $\R$ to $\R$.

We obtained a system which is not control-reducible but observation-reducible once.
Moreover, the corresponding pivot operator $[\A\ored{1}]$ is invertible, so all the pivot operators $[\A\ored{k}]$ are invertible.
This shows that \autoref{prop_findim_regnobeta} cannot hold in the infinite dimensional case.

As a result, the notion of \emph{index} is not appropriate for infinite dimensional systems.
The corresponding infinite dimensional notion is the data of all the combinations of observation and control reductions that lead to a reducible system.
In the finite dimensional case, if a system is regular it suffices to know the number of observation reductions (or control reductions), since both numbers are the same, and any combination will lead to the irreducible system.

This is no longer true in the infinite dimensional case.

Observe however that the number of combinations is fortunately limited by \autoref{thm_commutativity}, at least when all the reduced systems are normal.

\subsection{Saddle Point Problems}

Saddle Point Problems appear naturally in numerous applications (see, e.g., \cite{Brezzi}), and make a perfect example of operator pencil to study.
We will show that such systems are always normal, and we describe the effect of both reductions on them, thus confirming that reductions commute.
We also show that the invertibility of the pivot operators is precisely the inf-sup condition.

The ``stationary'' saddle point problem $\A u = f$ is studied in details in \autoref{sec_saddle_linear} 

\subsubsection{Definition}

Consider Hilbert spaces $X$ and $M$, and operator $A$ from $X$ to $X^*$, and $B$ from $X$ to $M^*$.
The corresponding saddle point problem is given by operators defined from 
\[
	\U :=  X \times M
\]
to 
\[
	\W := (X \times M)^* \equiv X^* \times M^*
	.
\]
The operators $\E$ and $\A$ are defined by
\begin{equation}
\label{eq_saddle_def}
	\E = \begin{bmatrix}
	R & 0 \\ 0 & 0
\end{bmatrix}
,
\qquad \A = \begin{bmatrix}
	A & B^* \\ B & 0
\end{bmatrix}
,
\end{equation}
where $B^*$ is defined as the \emph{transpose} of $B$,
%(terminology from \cite[15.1]{LaxFA}),
i.e., the operator defined from $M$ to $X^*$ defined by duality
\[
\pairing{B^* p}{u}_X := \pairing{B u}{p}_M \qquad \forall u,p \in X \times M
,
\]
and $R$ is the \emph{Riesz mapping} from the Hilbert space $X$ to its dual $X^*$.

\subsubsection{Normality}

We first show that systems stemming from saddle point problems are always normal.

\begin{proposition}
	\label{prop_saddle_normal}
	Systems of the form \eqref{eq_saddle_def} are normal.
\end{proposition}
\begin{pr}
We compute
\begin{equation}
	\label{eq_saddle_EU}
	\EU = \W\ored{1} = X^* \times 0 \equiv X^* 
	,
\end{equation}
and
\begin{equation}
	\label{eq_saddle_kE}
\ker\E = 0 \times M \equiv M
,
\qquad
\A\ker\E = B^* M
.
\end{equation}
As a result, since $\Ak \subset \EU$, the sytem is normal, as we observed in \autoref{rk_normal_no_alpha}.
	
\end{pr}

\subsubsection{Reduced Systems}

We now proceed to compute the remaining reduced subspace and corresponding pivot operators.

In addition to the spaces $\W\ored{1}$, $\kE$ and $\Ak$ computed in \eqref{eq_saddle_EU} and \eqref{eq_saddle_kE} we have
\[
\U\ored{1} = \ker B \times M
.
\]

For the observation reduced system we compute
\[
\ker\E\ored{1} = 0 \times M \equiv M
\qquad\text{and}\qquad
\A\ker\E\ored{1} \equiv B^* M
.
\]

Applying a control reduction to the last system yields the spaces
\[
\U\ocred = \U\ored{1}/\ker\E\ored{1} = \ker B
\qquad\text{and}\qquad
\W\ocred = W\ored{1}/A\ker\E\ored{1} = X^*/B^* M \equiv (B M)\orth
.
\]

The operator $[\A\cred{1}]$ is thus the operator $B$ with codomain $\A\ker\E$.

We also notice that
\[
\U/\U\ored{1} \equiv X/\ker B
,
\qquad \W/W\ored{1} \equiv M^*
.
\]
As a result, the operator $[\A\ored{1}]$ is the operator $B^*$ restricted on the subspace $X/\ker B$.

A control reduction of the original system yields the system given by the spaces
\[
\U\cred{1} \equiv X
\qquad\text{and}\qquad
\W\cred{1} = (X^* \times M^*)/B^* M \equiv (\ker B)^* \times M^*
.
\]

Applying an observation reduction yields
\[
\W\cored = \E\cred{1}\U\cred{1} \equiv (\ker B)^*
\qquad\text{and}\qquad
\U\cored = \{\, (u,p) : B^* p \in (\ker B)^* \,\} = \ker B
.
\]

\begin{figure}[H]
	\centering
\Uninemat{M \& M \& 0}{\ker B \times M \& X\times M \& M^* \& 0}{\ker B \& X \& M^* \& 0}
\Wninemat{0 \& 0}{M \& M \& 0}{X^* \& X^*\times M^* \& M^*}{X^*/B^* M \& X^*/B^* M \times M^* \& M^*}
\caption{
	The diagrams of \autoref{fig_commutativity} are shown here for the Saddle point problem, when the inf-sup condition (\ref{infsup}) is fulfilled.
All the rows and columns are exact.
}
\label{fig_saddle}
\end{figure}

\subsubsection{inf-sup Condition}

The \emph{inf-sup condition} (see, e.g., \cite[\S~4.1]{GiraultRaviart}, \cite[\S~II.2.3]{Brezzi}) for such a problem is the condition
\begin{equation}
\label{infsup}
\exists \beta>0 \qquad \inf_{\mu c\in M} \sup_{v\in X} \frac{\pairing{B^* \mu c}{v}}{\|\mu c\| \|v\|} \geq  \beta
.
\end{equation}

The inf-sup condition turns out to be exactly the condition of invertibility of the pivot operators.

\begin{proposition}
\label{prop_infsup}
The following statements are equivalent:
\begin{thmenumerate}
	\item the inf-sup condition \eqref{infsup} is fulfilled
	\item the operator $[\A\ored{1}]$ is invertible
	\item the operator $[\A\cred{1}]$ is invertible
\end{thmenumerate}
\end{proposition}
\begin{pr}
It is shown, for example, in \cite[Lemma 4.1]{GiraultRaviart}, that \eqref{infsup} is equivalent to the invertibility of $B^*$ from $M$ to its image $B^* M$ and that this image is closed.
This is in turn equivalent to the invertibility of $[\A\cred{1}]$.

Moreover, it is also shown that the inf-sup condition is equivalent to the fact that $B$ is invertible from $(\ker B)\orth$ to $M^*$, where $(\ker B)\orth$ is the subspace orthogonal to $\ker B$ with respect to the scalar product.
But that subspace $(\ker B)\orth$ is naturally isomorphic to $X/\ker B$, so the claim is proved.
\end{pr}

\subsection{``Index one'' Examples}

We turn to the study of some standard types of systems and show that they have ``control index one''.
Let us define that notion with the help of a straightforward proposition.

\begin{proposition}\label{prop_EUAK_empty}
	Consider a system $(\E,\A)$.
	The following statements are equivalent.
\begin{enumerate}
	\item The system $(\E_1,\A_1)$ is control-irreducible
		\item $\E\cred{1}$ is injective
			\item for all $u\in\U$, $\E u \in \AK \implies u \in \kE$
		\item $\E\U \cap \AK = 0$
\end{enumerate}

In that case we will say that the system has ``control index one''.
\end{proposition}

\subsubsection{Variational Systems}

In \cite[\S~4]{Tischendorf}, systems of a particular form are studied.
The domain $\U$ is a Hilbert space.
The codomain $\W$ is the dual space $\U^*$, i.e.,
\begin{equation*}
	\W := \U^*
	.
\end{equation*}

One is also given a ``pivot'' Hilbert space $\H$, such that $\U$ is densely included in $\H$.
With the identification $\H \equiv \H^*$ we obtain the Sobolev triple
\[
	\U \subset \H \equiv \H^* \subset \U^*
	.
\]

The operator $\E$ is of the form
\begin{equation}
	\label{eq_EDD}
	\E := \Do^* \Do
	,
\end{equation}
where the operator $\Do$ is defined from $\U$ to $\H$, and where we used the identification $\H \equiv \H^*$.

Finally, one assumes that
\begin{equation}
	\label{eq_Acoercive}
\A\quad\text{is coercive}
.
\end{equation}

We now show that such systems have ``control index one''.
\begin{proposition}
Consider a system $(\E,\A)$ such that \eqref{eq_EDD} and \eqref{eq_Acoercive} are fulfilled.
Assume moreover that $\A\kE$ is closed.
Then that system has control index one.
\end{proposition}
\begin{pr}
%Observe that a system $(\E,\A)$ is control-irreducible if and only if $\kE=0$.
%The claim will thus be proved if we show that $\ker\E\cred{1} = 0$.
\begin{enumerate}
		\item
Pick $u\in\U$ such that $\E u \in \AK$.
Since $\Ak$ is closed, there exists $k\in\kE$ such that $\E u = \A k$.
\item
We claim that $\pairing{\E u}{k} = 0$.
Indeed, by symmetry of $\E$, $\pairing{\Do^*\Do u }{k} = \ps{\Do u}{\Do k} = \pairing{\Do^* \Do k}{u} = 0$, since $k\in\kE$.
\item
This implies $\pairing{\A k}{k} = 0$ and thus $k=0$ by coercivity of $\A$.
As a result we obtain $\E u = 0$, so $u\in\kE$, which proves that $\ker\E\cred{1} = 0$, by \autoref{prop_EUAK_empty}.
\end{enumerate}
\end{pr}

Note that the property of $\A\kE$ to be closed is not restrictive.
It is in particular fulfilled whenever $[\A\cred{1}]$ is invertible, or in the finite dimensional case, and in general whenever $\rho (\E,\A)\neq \emptyset$ (see \autoref{lma_AkerE}).

\subsubsection{Inequality Constraints}
\label{sec_favini_indone}

In \cite[\S~2.4]{FaviniYagi}, the authors study systems for which $\W$ is a Hilbert space, and which fulfill the condition
\begin{equation}\label{eq_yagi}
	\exists \beta>0 \qquad \Re \ps{\E u}{\A u} \leq  \beta \| \E u \|^2 \qquad \forall u\in \U
.
\end{equation}
Let us show that under this condition the system studied must have ``control index one''.
This was pointed out in \cite[\S~4]{Tischendorf} but with different notations and definitions, so we show how this is true in our setting as well.

\begin{proposition}
	Under the assumption \eqref{eq_yagi}, the system $(\E,\A)$ has control index one.
\end{proposition}
\begin{pr}
By taking any element $u\in\U$ and $k\in\kE$ and assuming \eqref{eq_yagi} we obtain
\[
\Re \ps{\E(u+k)}{\A(u+k)} = \Re \ps{\E u}{\A k} + \Re \ps{\E u}{\A u} \leq  \beta \| \E u \|^2
.
\]
From this we see that we must have $\ps{\E u}{\A k} = 0$ for all $u\in\U$ and $k\in\kE$.
This implies in particular that
\[
\E\U \cap \A\kE = 0
,
\]
from which we conclude with \autoref{prop_EUAK_empty}.
\end{pr}

\section{Generalized Eigenvalue Problem}
\label{sec_generalized_spectrum}

We proceed to investigate problems of the kind
\[
	\lambda \E u + \A u = f
,
\]
for a given $f\in\W$.

Obviously, as a particular case, when $\lambda =0$, the problem reduces to
\[
	\A u = f
.
\]
We will discuss this type of problem in \autoref{sec_linprob}.

\subsection{Short Five Lemma}

We will use the Banach space version of a well known Lemma, generally used in homological algebra.
The proofs being obtained by diagram chasing, they are relatively easy to generalize to the Banach space case.
Some care is necessary, though, because operators on Banach spaces may have dense image without being surjective.

\begin{lemma}\label{lma_five_dec}
	Consider the notations and operators of \autoref{prop_subops}.
Then the following properties hold.
\begin{thmenumerate}
\item $\Sop$ injective $\implies$ $\Sop'$ injective
\item $\Sop$ surjective $\implies$ $[\Sop]$ surjective
\item $\Sop$ surjective and $[\Sop]$ injective $\implies$ $\Sop'$ surjective
\item $\Sop'$ surjective and $\Sop$ injective $\implies$ $[\Sop]$ injective
\item $[\Sop]$ surjective and $\Sop'$ surjective $\implies$ $\Sop$ surjective
\item $[\Sop]$ injective and $\Sop'$ injective $\implies$ $\Sop$ injective
\end{thmenumerate}
\end{lemma}

\begin{pr}
\begin{enumerate}
\item
The first claim is straightforward because $\ker\Sop' \subset \ker\Sop$.

\item
Similarly, consider $y+Y'\in Y/Y'$.
If $\Sop$ is surjective then there exists $x\in X$ such that $Sx = y$, so $[\Sop](x+X') = y + Y'$, and the second claim is proved.

\item
Consider $y\in Y'$.
Since $\Sop$ is injective, there exists $x\in X$ such that $\Sop x = y$.
Moreover, since $y + Y' = 0 + Y'$, by injectivity of $[\Sop]$, it follows that $x+X' = 0+X'$, i.e., that $x\in X'$, so $\Sop'$ is surjective.

\item
Assume that $\Sop(x+X') = 0+Y'$.
This is equivalent to $Sx \in Y'$.
Since $\Sop'$ is surjective, there exists $x_0\in X$ such that $\Sop x_0 = \Sop x$.
Since $\Sop$ is injective, $x_0 = x$, so $x+X' = x_0+X' = 0+X'$ and we conclude that $[\Sop]$ is injective.
\item
Consider $y\in Y$.
Since $[\Sop]$ is surjective, there exists $x\in X$ and $y'\in Y'$ such that $Sx = y + y'$.
Since $\Sop'$ is surjective, there exists $x'\in X'$ such that $Sx' = y'$, so $S (x-x') = y$ and the claim is proved.

\item
Consider $x\in X$ such that $\Sop x = 0$.
Since $[\Sop]$ is injective, this implies that $x\in X'$.
Using that $\Sop'$ is injective shows that $x=0$, which proves the claim.
\end{enumerate}
\end{pr}

The following corollary immediately follows from those results.

\newcommand*{\invert}{\quad\text{invertible}}

\begin{corollary}
Under the assumptions of \autoref{lma_five_dec} the following statements hold.
\begin{thmenumerate}
\item 
\[[\Sop] \invert \implies (\Sop \invert \iff \Sop' \invert)\]
\item 
\[\Sop' \invert \implies (\Sop \invert \iff [\Sop] \invert)\]
\end{thmenumerate}
\end{corollary}

\subsection{Resolvent Sets}

We are concerned with the generalized eigenvalue problem associated with the system $(\E,\A)$, i.e., for which $\lambda \in\CC$ the operator $\lambda \E + \A$ is invertible.

\begin{definition}
\label{def_resolvent_set}
The \emph{resolvent set} $\rho (\E,\A)$ of the system $(\E,\A)$ is defined by
\[
\rho (\E,\A) := \bigl\{\, \lambda \in\CC : \lambda \E + \A\quad\text{invertible} \,\bigr\}
.
\]
\end{definition}

\begin{lemma}\label{lma_AkerE}
If $\rho (\E,\A) \neq  \emptyset$ then $\A\kE$ is closed.
\end{lemma}
\begin{pr}
Observe that for any $\lambda  \in \CC$,
\[
	\A\kE = (\lambda \E + \A)\kE
	.
\]
Now if we choose $\lambda \in \rho (\E,\A)$, the operator $(\lambda \E+\A)$ has a continuous inverse, so it maps closed subspaces to closed subspaces.
The claim is now proved since $\kE$ is closed.
\end{pr}

\begin{theorem}
\label{thm_resolvent}
The following statements hold
\begin{enumerate}
	\item $[\A\ored{1}]$ invertible $\implies$ $\rho (\E,\A) = \rho (\E\ored{1},\A\ored{1})$
	\item $[\A\cred{1}]$ invertible $\implies$ $\rho (\E,\A) = \rho (\E\cred{1},\A\cred{1})$	
	\item $[\A\ored{1}]$ not invertible $\implies$ $\rho (\E,\A) = \emptyset$
	\item $[\A\cred{1}]$ not invertible $\implies$ $\rho (\E,\A) = \emptyset$
\end{enumerate}
\end{theorem}

\begin{pr}
The proof is a consequence of the following remarks
\begin{enumerate}
\item
For $\lambda \in\CC$, consider $\Sop_{\lambda } := \lambda \E + \A$.
We are going to use \autoref{lma_five_dec}.
As a subspace $\U' \subset \U$, we may either use $\U\ored{1}$ or $\kE$.
Using the notations of \autoref{lma_five_dec}, in the first case, $[\Sop_\lambda ] = [\A\ored{1}]$ and does not depend on $\lambda $, and in the second case, $\Sop'_{\lambda } = [\A\cred{1}]$ and does not depend on $\lambda $ either.
\item
Notice that $[\A\ored{1}]$ is always injective, so if it is not invertible, it is not surjective, and $\sigma (\E,\A)=\CC$.
\item
Notice that $[\A\cred{1}]$ has always a dense image, but needs not be surjective.
If $[\A\cred{1}]$ is not injective then $\rho (\E,\A) = \emptyset$.
\item The case of $[\A\cred{1}]$ not surjective is not covered by \autoref{lma_five_dec}.
We need to show that $[\A\cred{1}]$ not surjective implies that $\rho (\E,\A)=\emptyset$.
We conclude this by noticing that $[\A\cred{1}]$ not surjective is equivalent to $\A\kE$ not being closed, and using \autoref{lma_AkerE}.
\end{enumerate}
\end{pr}

\begin{remark}
\label{rk_empty_resolvent}
As we mentioned in \autoref{sec_findim_reg}, one of the obstacle in linking the invertibility of the pivot operators with the non-emptiness of the resolvent set is that the resolvent set of an irreducible system may be empty.

Consider for instance a Banach space $\U$ densely and continuously injected in a Banach space $\W$ but not closed in $\W$.
Choose $\E$ to be that injection and define $\A := 0$.
Such a system is irreducible and has an empty resolvent set, since $\E$ is not surjective.
\end{remark}

\section{Linear Problems}
\label{sec_linprob}

\subsection{General Result}

Let us briefly discuss problems of the kind
\[
	\A u = f
	,
\]
for an operator $\A$ defined from $\U$ to $\W$ and some element $f\in\W$.

\begin{corollary}
\label{cor_regular}
Pick an operator $\E$ with same domain and codomain as $\A$.
The system $(\overline{\E},\overline{\A})$ denotes the system $(\E,\A)$ after a finite number of reductions (observation or control).

Then, if any two of the following assertions hold, the third one holds as well.
\begin{itemize}
	\item $\overline{\A}$ is invertible
		\item $\A$ is invertible
		\item the pivot operators of the reductions leading to the system $(\overline{\E},\overline{\A})$ are all invertible
\end{itemize}
\end{corollary}
\begin{pr}
	The proof is a straightforward consequence of \autoref{thm_resolvent}, since $\A$ is invertible if and only if $0\in \rho (\E,\A)$.
\end{pr}

Note that there is no need to make any assumptions on the system $(\overline{\E},\overline{\A})$, neither on the normality of any of the reduced systems.

The power of \autoref{cor_regular} depends on a judicious choice of the auxiliary operator $\E$.
Let us discuss some extreme choices.

The first extreme is to choose $\E = 0$, which yields one reduction steps (either observation or control) and the pivot operator in either case is $\A$ itself.
In that case, \autoref{cor_regular} reduces to a tautology, since $\overline{\A}$ is empty, i.e. the zero operator from the zero-dimensional vector space to itself, thus invertible.

The other extreme choice is, if possible, to choose $\E$ injective and with a dense image.
In that case, the system is irreducible, and the pivot operators are empty, thus invertible, and $\overline{\A} = \A$.

\subsection{Saddle Point Problem}
\label{sec_saddle_linear}

Now we obtain the standard result of Saddle point problems using \autoref{thm_commutativity}.
Define the injection $\ii$ by
\[
\ii:\quad \ker B \rightarrow X
.
\]

We obtain the following standard result (\cite[Theorem~4.1]{GiraultRaviart}).
\begin{proposition}
The operator $\A$ is invertible if and only if $\ii^* A \ii$ is invertible and the inf-sup condition \eqref{infsup} is fulfilled.
\end{proposition}
\begin{pr}
The operator $\ii^* A \ii $ is none other than the operator $\A\ored{1}\cred{1}$.
If the inf-sup condition is fulfilled, then, using \autoref{thm_commutativity}, we obtain that $[\A\ocred]$ is invertible.
We may then use \autoref{cor_regular} to obtain that $\A$ is invertible if and only if $\A\ocred$ is invertible, i.e., if and only if $\ii^* A \ii$ is invertible.

On the other hand, if $\ii^* A \ii$ is invertible and $\A$ is invertible, then by \autoref{cor_regular}, we obtain that both $[\A\cred{1}]$ must be invertible. This is equivalent to the inf-sup condition by \autoref{prop_infsup}.
\end{pr}

\subsection{Formulations of the Poisson Problem}

Let us consider the case of the Poisson problem on bounded open set $\Omega \subset \R^d$.
We will show that the ``mixed'' formulation of the Poisson problem corresponds to a special choice of an auxiliary operator $\E$.

The Poisson problem may be interpreted as the stationary version of the heat equation.
We thus express the heat equation as a first order evolution equation.

With the convention that all the spaces are defined on the domain $\Omega $, the spaces $\U$ and $\W$ are defined by
\[
	\U := \H_0^1 \times \Leb2
	,
	\qquad
	\W := \U^* \equiv  \H\inv \times \Leb2
	.
\]

The operators $\E$ and $\A$ are then defined as
\[
	\E = \begin{bmatrix}
		\ii^* \ii & 0 \\ 0 & 0
	\end{bmatrix}
	,
	\qquad
	\A = \begin{bmatrix}
		0 & \div \\ \grad & R
	\end{bmatrix}
	,
\]
where $R$ is the Riesz mapping from $\Leb2$ to its dual.

The system $(\E,\A)$ may be control-reduced, or observation-reduced once.
In both cases, the reduced operators $\A\ored{1}$ and $\A\cred{1}$ may be interpreted as the Laplace operator from $\H_0^1$ to $\H\inv$.

Since this brings nothing new, let us change the setting.
Setup instead
\[
	\U := \Leb2 \times \H(\div)
	,
	\qquad
	\W := \U^*
.
\]
The operators are now given by
\[
	\E := \begin{bmatrix}
		0 & 0 \\ 0 & R
	\end{bmatrix}
	,
	\qquad
	\A := \begin{bmatrix}
		0 & \div \\ \grad & \ii^* \ii
	\end{bmatrix}
	,
\]
where $R$ is now the Riesz mapping from $\H(\div)$ in its dual, and $\ii$ is the injection from $\H(\div)$ to $\Leb2$.

In that case, the system is a saddle point problem (see \cite{Brezzi}).

\section{Conclusion and Outlook}

We have studied in detail the effect of observation and control reduction on system of operators on Banach spaces.
The main result is \autoref{thm_commutativity}, according to which those reduction commute under the normality assumptions of \autoref{def_normality}.

What is the structure of existing implicit differential equations from the reduction point of view?
We studied the saddle point problem, the multiplication operator in \autoref{sec_multiplication_operator}, as well as various systems appearing in \cite{FaviniYagi} and \cite{Tischendorf}.
However, there are many other systems of interest to be studied.
Let us mention for instance the linearized elastodynamics in \cite{Simeon}, the Dirac equation in the nonrelativistic limit in \cite[\S~3]{ThallerCauchy}, linear PDEs as studied in \cite{Campbell} and \cite{Seiler}.
In the latter cases, it would be interesting to compare the reduction structure that we obtain to the index concept developed, in particular in \cite{Campbell}.

This brings us to an essential question concerning operator pencil: is there an equivalent of the Kronecker decomposition theorem?
What we did in this work was to inspect two consequences of the Kronecker theorem, and examine their validity in the infinite dimensional case.
According to \autoref{thm_commutativity}, the commutativity of reduction, which is a consequence of the Kronecker decomposition in the finite dimensional case, is still true in the infinite dimensional case, at least under some conditions.
The counter example of \autoref{sec_mult_reg_ored} shows however that some other consequences of the Kronecker decomposition theorem, namely \autoref{prop_findim_regnobeta}, which is essential to define the notion of \emph{index}, do not hold anymore in the infinite dimensional case.

The question remains of which other structures from the finite dimensional case are preserved and much more remains to do in that respect.
%For instance, is there an equivalent of the \emph{constraint defects} which were defined in \autoref{sec_findim} (see \cite{Verdier})?

\subsubsection*{Acknowledgements}
	
I would like to acknowledge the support of the \href{http://wiki.math.ntnu.no/genuin}{GeNuIn Project}, funded by the \href{http://www.forskningsradet.no/}{Research Council of Norway}, and that of its supervisor, \href{http://www.math.ntnu.no/~elenac/}{Elena Celledoni}.
I also thank the anonymous referees for their helpful remarks.

\bibliographystyle{abbrvurl}
\bibliography{linpdae}

\end{document}